\documentclass{amsart}
\usepackage{latexsym,amsfonts,amsmath,amssymb}

\newcommand{\IR}{{\mathbb{R}}}
\newcommand{\ID}{{\mathbb{D}}}

\newcommand{\IC}{{\mathbb{C}}}

\newcommand{\C}{{\mathcal{C}}}
\newcommand{\E}{{\mathcal{E}}}

\newcommand{\LL}{{\mathcal{L}}}
\newcommand{\M}{{\mathcal{M}}}

\newcommand{\al}{{\alpha}}
\newcommand{\alt}{{\tilde\alpha}}

\newcommand{\Tr}{{\rm{Tr}}}

\newcounter{smalllist}

\newtheorem{theorem}{Theorem}

\newtheorem{lemma}{Lemma}[section]

\newtheorem{coro}[lemma]{Corollary}
\theoremstyle{definition}

\newtheorem{remark}[lemma]{Remark}

\DeclareMathOperator{\diag}{diag} 

\let\llldots=\ldots
\def\ldots{\llldots{}}

\numberwithin{equation}{section}

\begin{document}

\title{A note on circular trace formulae}
\author{Irina Nenciu}
\address{Irina Nenciu\\
         School of Mathematics\\
         Institute for Advanced Study\\
         Princeton, NJ 08540}
\email{nenciu@ias.edu}
\thanks{The author wishes to thank Barry Simon, for his encouragement and useful suggestions.
This work was partly supported by NSF grant DMS-0111298.}

\begin{abstract}
We find a finite CMV matrix whose eigenvalues coincide with the
Dirichlet data of a circular periodic problem. As a consequence,
we obtain circular analogues of the classical trace formulae for periodic
Jacobi matrices.
\end{abstract}

\maketitle

\section{Introduction and background}

CMV matrices came to the fore in the theory of orthogonal polynomials on the unit circle (OPUC)
within the past five years as a unitary analogue of Jacobi matrices.
This analogy manifests itself in many fields, from orthogonal polynomials and numerical
analysis, to random matrices and integrable systems. For more details, we direct the reader to the survey papers
\cite{Nen} and \cite{Simon3}, and the references therein. In this note, we present a circle analogue
of the well-known trace formulae for periodic Jacobi matrices (see \cite{vanM}), thus answering a question left open by
Simon, \cite{Simon2}.

We begin with a brief overview of the theory of OPUC
and CMV matrices; for a more complete description of what follows, the reader
may turn to \cite{Simon1} and \cite{Simon2}. In this note, we will use the definitions and notation
of \cite{Simon1} and \cite{Simon2}; the reader familiar with the material can skip ahead to
Section~2.

Given a probability measure $d\mu$ on $S^1$, the
unit circle in $\IC$, we can construct an orthonormal system of
polynomials, $\phi_k$, by applying the Gram--Schmidt procedure to
$\{1,z,\ldots\}$.  These obey a recurrence relation; however, to
simplify the formulae, we will present the relation for the monic
orthogonal polynomials $\Phi_k$:
\begin{align}
\Phi_{k+1}(z)   &= z\Phi_k(z)   - \bar\alpha_k \Phi_k^*(z).
\label{PhiRec}
\end{align}
Here, $\alpha_k\in\IC$ are recurrence coefficients, which are called
Verblunsky coefficients, and $\Phi_k^*$ denotes the reversed
polynomial:
\begin{equation}\label{E:rev}
\Phi_k(z) = \sum_{l=0}^k c_l z^l \quad \Rightarrow \quad \Phi_k^*(z)
= \sum_{l=0}^k \bar{c}_{k-l} z^l=z^k\overline{\Phi_k\left(\frac{1}{\bar z}\right)}.
\end{equation}

When $d\mu$ is supported at exactly $n$ points, $\alpha_k\in\ID=\{z\in\IC\,|\,|z|<1\}$
for $0\leq k\leq {n-2}$ while $\alpha_{n-1}$ is a unimodular complex
number. If the support of $d\mu$ is infinite, then the recurrence formula \eqref{PhiRec}
will produce an infinite sequence $\{\al_k\}_{k\geq0}$ of Verblunsky coefficients, all of whom are inside the unit disc $\ID$.
In both of these cases, there is a 1-1 correspondence between the measure and the sequence of coefficients.

The situation described above parallels the one on the real line $\IR$.
Given a probability measure $d\nu$ supported on $\IR$, we can apply the Gram--Schmidt
procedure to $\{1,x,x^2,\ldots\}$ and so obtain an
orthonormal basis for $L^2(d\nu)$ consisting of
polynomials, $\{p_j\}_{j\geq0}$, with positive leading
coefficient. In this basis, the linear transformation $f(x)\mapsto
xf(x)$ is represented by a Jacobi matrix: a tri-diagonal matrix
\begin{equation}\label{Jmat}
J=\begin{bmatrix}
 b_1  &  a_1 &        & \\
 a_1  &  b_2 & a_2 &  \\
      &a_2& b_3 & \ddots\\
      &      & \ddots&  \ddots
\end{bmatrix}
\end{equation}
with $a_j>0$, $b_j\in\IR$.
An equivalent statement
is that the orthonormal polynomials obey a three-term recurrence relation:
$$
xp_{j}(x) = a_{j+1} p_{j+1}(x) + b_{j+1}p_j(x) + a_{j} p_{j-1}(x)
$$
where $a_{0}=0$. As in the circle case, if the support of $d\nu$ consists
of exactly $n$ points, then $a_n=0$, and $J$ becomes an $n\times n$
symmetric matrix.
We have just shown how measures on $\IR$ lead to Jacobi matrices;
in fact, there is a one-to-one correspondence between them.  Given a
Jacobi matrix, $J$, let $d\nu$ be the spectral measure associated to
$J$ and the vector $e_1=[1,0,\ldots,0]^{T}$.  Then $J$ represents
$x\mapsto xf(x)$ in the basis of orthonormal polynomials associated
to $d\nu$.

From the discussion of Jacobi matrices, it would be natural to consider the matrix representation
of $f(z)\mapsto zf(z)$ in $L^2(d\mu)$ with respect to an appropriate basis.
Cantero, Moral, and Vel\'azquez, \cite{CMV}, had the simple and ingenious idea of applying the Gram--Schmidt
procedure to $\{1,z,z^{-1},z^2,z^{-2},\ldots\}$ rather than $\{1,z,\ldots\}$.
In the resulting basis, the map $f(z)\mapsto zf(z)$ is represented in an especially simple form:
Given the coefficients $\alpha_0,\al_1,\ldots$ in $\ID$, let $\rho_k=\sqrt{1-|\alpha_k|^2}$, and
define $2\times 2$ matrices
$$
\Theta_k = \begin{bmatrix} \bar\alpha_k & \rho_k \\ \rho_k & -\alpha_k
\end{bmatrix}
$$
for $k\geq0$, while $\Theta_{-1}=[1]$ is a $1 \times1$ matrix. From these,
form the  block-diagonal matrices
$$
\mathcal{L}=\diag\bigl(\Theta_0   ,\Theta_2,\Theta_4,\ldots\bigr)
\quad\text{and}\quad
\mathcal{M}=\diag\bigl(\Theta_{-1},\Theta_1,\Theta_3,\ldots\bigr).
$$
Then the operator of multiplication by $z$ in $L^2(d\mu)$ can be represented
by the associated \emph{CMV matrix}, $\C=\mathcal{LM}$. If the support of
$d\mu$ consists of exactly $n$ points, then $\rho_{n-1}=0$, and so $\Theta_{n-1}=[\bar\al_{n-1}]$
becomes a $1\times1$ unitary matrix, and the associated (finite) CMV matrix will be an $n\times n$
unitary matrix.

\begin{remark}\label{R:1}
Let us make two observations. First of all, if we consider a general (infinite)
CMV matrix, and set one of the Verblunsky coefficients to be on the unit circle $S^1$,
then the matrix decomposes as a direct sum. So a finite CMV matrix, where
$\al_{-1},\al_{n-1}\in S^1$, is in fact a unitary operator with Dirichlet boundary conditions.
The second observation is historical: Within the orthogonal polynomials
community, CMV matrices were indeed introduced by Cantero, Moral, and Vel\'azquez,
and hence the name of these special unitary matrices. But as it turns out, these matrices have been known in the numerical analysis community
for some 20 years; for more details, see, for example, \cite{Simon3} and \cite{Wat}.
\end{remark}

We are interested in sequences of Verblunsky coefficients
$\{\al_j\}_{j\geq0}$ which are periodic with period $p$:
$$
\al_{j+p}=\al_j\qquad\text{for all}\quad j\geq0.
$$
These are completely described by their first $p$ terms, so from now on,
whenever we talk about periodic Verblunsky coefficients, we will think
of finite sets $\{\al_j\}_{j=0}^{p-1}\in\ID^p$; for simplicity, we will assume throughout
this paper that the period $p$ is even.
Let us first observe that in the periodic case, besides the usual CMV matrix,
one can also define a so-called extended CMV matrix, that we
shall denote by $\E$. Indeed, starting with $\al_0,\ldots,\al_{p-1}\in\ID$, we
can define by periodicity a two-sided infinite sequence of coefficients. The extended CMV matrix is
\begin{equation*}\label{DefE}
\mathcal{E}=\mathcal{\tilde{L}\tilde{M}},
\end{equation*}
where
\begin{equation}\label{LMtilda}
\mathcal{\tilde L}=\bigoplus_{j\,\,\text{even}} \Theta_j \quad\quad
\text{and} \quad\quad \mathcal{\tilde M}=\bigoplus_{j\,\,\text{odd}}
\Theta_j,
\end{equation}
with $\Theta_j$ defined on $l^2(\mathbb{Z})$ by
$$
\Theta_j = \begin{bmatrix} \bar\al_j & \rho_j \\ \rho_j & -\al_j
\end{bmatrix}
$$
on the span of $\delta_j$ and $\delta_{j+1}$, and identically 0
otherwise.

Note that $\E$ acts boundedly on the space of bounded sequences $l^{\infty}$.
Moreover, if $S$ is the $p$-shift $(Su)_m=u_{m+p},\quad\text{for}\,\, u\in l^{\infty}$,
then, by periodicity of the $\al$'s, we see that $S\E=\E S$.
In particular, if $\beta\in S^1$ and we consider
$$
X_{\beta}=\{u\in l^{\infty}\,|\,Su=\beta u\},
$$
then $\E$ takes $X_{\beta}$ to itself:
$\E(X_{\beta})\subset X_{\beta}.$

We can therefore define
$$
\E(\beta)=\E\upharpoonright X_{\beta}.
$$
If we represent this operator in the natural basis in $X_{\beta}$, we obtain
the Floquet CMV matrix:
$$
\E(\beta)=\LL_p\M_p(\beta),
$$
with
$$
\LL_p=
\left(%
\begin{array}{ccccc}
  \Theta_0 &  &  &  &  \\
   & \ddots &  &  &  \\
   &  & \ddots &  & \\
   & & & \ddots&  \\
   &  &  &  & \Theta_{p-2} \\
\end{array}%
\right)
$$
and
$$
\M_p(\beta)=
\left(%
\begin{array}{ccccc}
  -\al_{p-1} &  &  &  & \rho_{p-1}\beta^{-1} \\
   & \Theta_1 &  &  &  \\
   &  & \ddots &  & \\
   & & & \Theta_{p-3}&  \\
   \rho_{p-1}\beta&  &  &  & \bar\al_{p-1} \\
\end{array}%
\right).
$$
Recall that we assumed $p$ is even.

Let $d\mu$ be the measure associated to the 1-sided sequence of periodic Verblunsky coefficients
defined by $\al_0,\ldots,\al_{p-1}$.
Then $d\mu$ is supported on the union of $p$ bands, $B_1,\dots,B_p$, on which the
measure is absolutely continuous, with at most one pure point between any two neighboring bands.
The readers familiar with the real line case will recognize this picture. The bands are usually described
in terms of the discriminant associated to the periodic problem; here we take a different point of
view, and claim that $\bigcup_{j=1}^p B_j=\bigcup_{\beta\in S^1} \text{spec}\bigl(\E(\beta)\bigr)$ (see
Theorem~11.1.1 and formula (11.2.17) in \cite{Simon2}).
Furthermore, the edges of the bands are given by the eigenvalues of $\E(\pm1)$. Let
$\text{spec}\bigl(\E(\pm1)\bigr)=\{z_1^\pm,\ldots,z_p^\pm\}\subset S^1$. If one travels on the unit circle
in a counterclockwise direction, and starting from an appropriate point, one encounters
these eigenvalues in the order $z_1^+,z_1^-,z_2^-,z_2^+,z_3^+,z_3^-,\ldots$, and the bands are
given by the (small) arcs $[z_{2j-1}^+,z_{2j-1}^-]$ and $[z_{2j}^-,z_{2j}^+]$. The possible
positions of the pure points are given by the zeroes $z_1,\ldots,z_p$ of the polynomial $\Phi_p-\Phi_p^*$,
where we recall that the $\Phi_k$'s, $k\geq 0$, are the monic orthogonal polynomials associated
to the measure $d\mu$, and $\Phi_k^*$
denotes the reversed polynomial (see \eqref{E:rev}). Let $\{\zeta_k\}_{k=1}^{2p}=\{z_j^\pm\}_{j=1}^p$,
so that the gaps (arcs that separate the bands) are given by $G_j=(\zeta_{2j-1},\zeta_{2j})$,
with $z_j\in G_j$, $1\leq j\leq p$.

The analogous points that give the putative positions of the pure points for periodic Jacobi matrices
are obtained as the eigenvalues of a Jacobi matrix with Dirichlet boundary conditions, and hence,
together with the information of which of then are indeed pure points, bear the name of Dirichlet data.
When studying the case of periodic Verblunsky coefficients, Simon \cite{Simon2} calls
the pairs $(z_j,\sigma_j)$ Dirichlet data, where the $\sigma_j$'s are $\pm1$,
depending on whether the corresponding $z_j$ is or is not a pure point of the measure $d\mu$.
But in this case the name was given just by analogy to the self-adjoint case, and the question
of the existence of a unitary operator with Dirichlet boundary conditions having the
$z_j$'s as eigenvalues was left open (see the Remarks and Historical Notes at the end of Section~11.3 of \cite{Simon2}).
This is exactly the question we answer in the second part of this note.

\section{Dirichlet eigenvalues and trace formulae}\label{S:2}

As explained in the Introduction, we focus here on the case in which the sequence of Verblunsky coefficients is periodic,
with period $p$: $\al_{j+p}=\al_j$ for all $j\geq 0$. For simplicity, we assume that $p$ is even.
The support of the associated
measure $d\mu$ will consist of $p$ bands, separated by gaps, and at most one pure point in each gap.
These points are a subset of the set $\{z_1,\ldots,z_p\}$ of zeroes of the polynomial
$\Phi_p-\Phi_p^*$. When studying the properties of periodic Verblunsky coefficients,
Simon \cite{Simon2} calls the $z_j$'s, together with the information of which are the pure points
of the measure, Dirichlet data. The name is given by analogy with the Jacobi case, where the positions of the pure points are actually
eigenvalues of an operator with Dirichlet boundary conditions.

The main result of this note is the following

\begin{theorem}\label{T:Dir}
Let $\{\al_j\}_{j=0}^{p-1}\in \ID^p$ determine a sequence of $p$-periodic Verblunsky coefficients,
and $z_1,\dots,z_p\in S^1$ be the positions of the associated Dirichlet data. Define a finite
sequence of Verblunsky coefficients $\{\alt_j\}_{j=0}^{p-1}\in\ID^{p-1}\times S^1$ by
$$
\alt_j=\al_j,\quad 0\leq j\leq p-2,\qquad \alt_{p-1}=\frac{1+\al_{p-1}}{1+\bar\al_{p-1}},
$$
and let $\C$ be the associated finite CMV matrix. Then the spectral measure associated to
$\C$ and the vector $e_1=[1\,0\dots 0]^T$ is $d\tilde\mu=\sum_{j=1}^p w_j \delta_{z_j}$, where
\begin{equation}\label{E:masses}
w_j=\left[\left.\frac{d}{dr}\right|_{r=1}|\phi_p(rz_j)|^2-p|\phi_p(z_j)|^2\right]^{-1}
\end{equation}
and $\phi_p$ is the $p^\text{th}$ orthonormal polynomial for the periodic problem.
\end{theorem}

\begin{remark}
In particular, this theorem tells us that $z_1,\ldots,z_p$ are the eigenvalues
of the finite CMV matrix $\C$, thus justifying the name of Dirichlet data (see also Remark~\ref{R:1}).
But, as is usually the case with circle analogues of real line quantities and phenomena,
the algebra becomes more involved: indeed, while only one of the coefficients needs to be changed,
same as in the real line case, this affects eight entries of the Floquet CMV matrix $\E(\beta)$, two of which are on the diagonal,
versus only two entries of the analogous Jacobi matrix. The effects of
this phenomenon will be very transparent in the statement of Corollary~\ref{C:Tr}.
\end{remark}

\begin{proof}
The fact that the $z_j$'s are the eigenvalues of the finite CMV matrix $\C$ follows immediately from the observation that
\begin{equation}\label{E:polys1}
\Phi_{p}(z)-\Phi_{p}^*(z)
=(1+\al_{p-1})\left[z\Phi_{p-1}(z)
-\frac{1+\bar\al_{p-1}}{1+\al_{p-1}}\Phi_{p-1}^*(z)\right]
\end{equation}
where we use the recurrence relation \eqref{PhiRec} for the periodic problem. Since
$\left|\frac{1+\al_{p-1}}{1+\bar\al_{p-1}}\right|=1$, the $\tilde\al_j$'s defined in the statement
can indeed be considered as the Verblunsky coefficients of a finite measure.
Let $\tilde\Phi_j$, $0\leq j\leq p$ be the associated
monic orthogonal polynomials. Since an orthogonal polynomial of degree $k$ depends only on the coefficients
of index 0 through $k-1$, and $\alt_j=\al_j$ for $j=0,\ldots,p-2$, we see that
$$
\tilde\Phi_k(z)=\Phi_k(z),\qquad\text{for all}\quad0\leq k\leq p-1.
$$
So the last recurrence relation for the finite problem is
\begin{align*}
\tilde\Phi_p(z) &=z\Phi_{p-1}(z)-\overline{\alt_{p-1}}\Phi^*_{p-1}(z)\\
                &=z\Phi_{p-1}(z)-\frac{1+\bar\al_{p-1}}{1+\al_{p-1}}\Phi_{p-1}^*(z),
\end{align*}
and hence \eqref{E:polys1} becomes
\begin{equation}\label{E:polys2}
\Phi_{p}(z)-\Phi_{p}^*(z)=(1+\al_{p-1})\tilde\Phi_p(z).
\end{equation}
But $\tilde\Phi_p(z)=\det(z-\C)$, and therefore the $z_j$'s are the eigenvalues of $\C$, or, equivalently,
$\{z_1,\dots,z_p\}$ is the support of the spectral measure of $\C$ and $e_1$.

Theorem~11.5.8 of \cite{Simon2} tells us that, if we define the $w_j$'s as in the statement of the theorem,
for $j=1,\ldots,p$, then $\phi_0,\ldots,\phi_{p-1}$, the first $p$ orthonormal polynomials
of the periodic problem, are orthonormal with respect to
the (trivial) measure $d\tilde\mu=\sum_{j=1}^p w_j\delta_{z_j}$. If we set $\beta_0,\ldots,\beta_{p-2}\in\ID$
and $\beta_{p-1}\in S^1$ to be the Verblunsky coefficients of $d\tilde\mu$, it follows from the recurrence
relation for $d\tilde\mu$ that
$\beta_0=\alpha_0,\ldots,\beta_{p-2}=\alpha_{p-2}$. Moreover, by taking the determinant of
the associated finite CMV matrix, we find that $-\bar\beta_{p-1}=\prod_{j=1}^p z_j$ (see, for example, \cite{KN}).
But by the same argument applied to the matrix $\C$, we obtain that
$$
\frac{1+\alpha_{p-1}}{1+\bar\alpha_{p-1}}=-\prod_{j=1}^p \bar z_j=\beta_{p-1}.
$$
In other words, $d\tilde\mu$ is exactly the measure associated with the finite CMV matrix $\C$, which concludes the proof.
\end{proof}

As a simple consequence of Theorem~\ref{T:Dir}, one obtains the following trace formulae:
\begin{coro}[Trace formulae]\label{C:Tr}
Let $\{\al_j\}_{j=0}^{p-1}\in \ID^p$ determine a sequence of $p$-periodic Verblunsky coefficients,
$\{\zeta_k\}_{k=1}^{2p}=\{z_j^\pm\}_{j=1}^p$ be the endpoints of the gaps,
and $z_1,\dots,z_p\in S^1$ be the positions of the associated Dirichlet data. Then the following
relations hold:
\begin{equation}\label{E:det1}
\prod_{j=1}^p z_j=-\frac{1+\al_{p-1}}{1+\bar\al_{p-1}}
\end{equation}
and
\begin{equation}\label{E:Tr1}
\sum_{j=1}^p \Bigl(\frac{\zeta_{2j-1}+\zeta_{2j}}{2}-z_j\Bigr)
=-\bar\al_0 (1+\al_{p-1}) +\al_{p-2}\frac{1-|\al_{p-1}|^2}{1+\al_{p-1}}.
\end{equation}
\end{coro}

\begin{proof}
The first equation follows immediately by taking the determinant of the finite CMV matrix $\C$. This relation appears already
in \cite{Simon2} (see formula (11.3.24)), but with a different proof. Recalling that
$\prod z_j^+=\prod z_j^-=1$, since the $z_j^\pm$'s are the eigenvalues of unitary matrices, we can rewrite \eqref{E:det1} as
\begin{equation}\label{E:det2}
\prod_{j=1}^p \bigl(\zeta_{2j-1}\zeta_{2j}z_j^{-2}\bigr)=\left(\frac{1+\bar\al_{p-1}}{1+\al_{p-1}}\right)^2.
\end{equation}
Note that, if we translate this relation in terms of the arguments of the $\zeta_j$'s and $z_j$'s, we
obtain a close analogue of the first trace formula for periodic Jacobi matrices.

The second equation represents $\frac12\Tr\bigl(\E(+1)+\E(-1)-2\C\bigr)$ written in terms of the corresponding eigenvalues
on the left-hand side, and of the Verblunsky coefficients on the right-hand side.
\end{proof}

One can of course consider traces of higher powers, but the formulae quickly become quite intractable. As an example,
we write the next trace formula, for $\frac12 \Tr(\E(+1)^2+\E(-1)^2-2\C^2)$:
\begin{align*}
\sum_{j=1}^p \Bigl(\frac{\zeta_{2j-1}^2+\zeta_{2j}^2}{2}-z_j^2\Bigr)
&=\al_{p-2}^2\left(\bar\al_{p-1}^2+\Bigl(\frac{1+\bar\al_{p-1}}{1+\al_{p-1}}\Bigr)^2\right)-\bar\al_0^2(1-\al_{p-1}^2)\\
&\quad +2\al_{p-3}\rho_{p-2}^2\frac{\rho_{p-1}^2}{1+\al_{p-1}}-2\rho_0^2\bar\al_1(1+\al_{p-1})\\
&\quad -2\al_{p-2}\rho_{p-1}^2\bar\al_0
\end{align*}
Recall that, in the real line case, one can use the trace formulae to explicitly
write down the recurrence coefficients, $a_k$ and $b_k$, in terms of the positions
of the bands and the Dirichlet data (see, for example, \cite{vanM}). On the circle,
we are unable to do so, mainly because the algebra is much more involved (see also the comments
in Remark~2.1). So it would be interesting to find a direct way of expressing the $\al_k$'s
in terms of the $\zeta_j$'s and $z_j$'s. At the same time, one should bear in mind the
fact that trace formulae are not necessary to develop the theory of
periodic Jacobi or CMV matrices; for more comments, see Chapter 11 in \cite{Simon2}.

As a final remark, we briefly present the analogous results for Aleksandrov measures.
Indeed, let $\lambda\in S^1$ be a complex number on the unit circle. One can then consider
the problem defined by the rotated periodic Verblunsky coefficients $\{\al_{\lambda,j}=\lambda\alpha_j\}_{j=0}^{p-1}\in\ID^p$,
the associated measure $d\mu_\lambda$ (called Aleksandrov measure), and
the monic orthogonal polynomials $\Phi_{\lambda,k}$.
It is easy to show that, while the bands of $d\mu_\lambda$ are $\lambda$-independent,
the positions $z_{\lambda,j}$ of the Dirichlet data (or, equivalently, the zeroes of
the associated polynomial $\Phi_{\lambda,p}-\Phi^*_{\lambda,p}$) change with $\lambda$. The same arguments as above
yield the fact that
the zeroes of the polynomial $\Phi_{\lambda,p}-\Phi_{\lambda,p}^*$
coincide with the eigenvalues of the finite CMV matrix $\C_\lambda$
defined by the Verblunsky coefficients
$\tilde\al_{\lambda,0}=\lambda\al_0,\ldots,\tilde\al_{\lambda,p-2}\lambda\al_{p-2}\in\ID$ and
$\tilde\al_{\lambda,p-1}=\frac{1+\lambda\al_{p-1}}{1+\bar\lambda\bar\al_{p-1}}\in S^1$. In particular,
$\Phi_{\lambda,p}-\Phi_{\lambda,p}^*$ has $p$ simple zeroes,
$z_{\lambda,1},\ldots, z_{\lambda,p}$, all of whom lie on the
unit circle.
The spectral measure associated to the finite CMV matrix $\C_\lambda$ is
$$\tilde\mu_{(\lambda)}=\sum_{j=1}^p w_{\lambda,j}\delta_{z_{\lambda,j}},
$$
where
$$
w_{\lambda,j}=\left[\left.\frac{d}{dr}\right|_{r=1}|\phi_{\lambda,p}(rz_{\lambda,j})|^2-p|\phi_{\lambda,p}(z_{\lambda,j})|^2\right]^{-1}
$$
Note that $d\tilde\mu_{(\lambda)}$ is NOT the Aleksandrov measure $d\tilde\mu_\lambda$ of $d\tilde\mu$, since
$\tilde\al_{\lambda,j}=\lambda\tilde\al_j$ for $j=0,\ldots,p-2$, but
$$
\frac{1+\lambda\al_{p-1}}{1+\bar\lambda\bar\al_{p-1}}\neq \lambda\cdot\frac{1+\al_{p-1}}{1+\bar\al_{p-1}}.
$$
Taking the appropriate traces leads to:
\begin{coro}[Trace formulae for Aleksandrov measures]\label{C:TrL}
Let $\{\al_j\}_{j=0}^{p-1}\in \ID^p$ determine a sequence of $p$-periodic Verblunsky coefficients,
$\{\zeta_k\}_{k=1}^{2p}=\{z_j^\pm\}_{j=1}^p$ be the endpoints of the gaps,
and $z_{\lambda,1},\dots,z_{\lambda,p}\in S^1$ be the positions of the Dirichlet data for the
associated Aleksandrov measure. Then the following
relations hold:
\begin{equation}\label{E:det1L}
\prod_{j=1}^p z_{\lambda,j}=-\frac{1+\lambda\al_{p-1}}{1+\bar\lambda\bar\al_{p-1}}
\end{equation}
and
\begin{equation}\label{E:Tr1L}
\sum_{j=1}^p \Bigl(\frac{\zeta_{2j-1}+\zeta_{2j}}{2}-z_{\lambda,j}\Bigr)
=-\bar\lambda\bar\al_0 (1+\lambda\al_{p-1}) +\lambda\al_{p-2}\frac{1-|\al_{p-1}|^2}{1+\lambda\al_{p-1}}.
\end{equation}
\end{coro}


\end{document}